\documentclass[12pt,letterpaper]{article}
\usepackage{mathrsfs}
\usepackage{indentfirst}
\usepackage{makeidx}
\usepackage{amsmath,amssymb,amsfonts,amsthm,graphics,graphicx,pifont}
\usepackage{makeidx}
\usepackage{float}
\usepackage{color}
\usepackage[font=small,labelfont=bf,labelsep=none]{caption}
\usepackage{epstopdf}

\newtheorem{thm}{Theorem}[section]

\newtheorem{lem}[thm]{Lemma}
\newtheorem{cor}[thm]{Corollary}
\newtheorem{rem}{Remark}[section]

\theoremstyle{definition}

\def\-{\mbox{--}}

\def\pf{\noindent {\it Proof.} }

\begin{document}
\title{Conflict-free (vertex)-connection numbers of graphs with small diameters\footnote{Supported by NSFC No.11871034, 11531011 and NSFQH No.2017-ZJ-790.}}
\author{
\small Xueliang Li$^{1,2}$, Xiaoyu Zhu$^1$\\
\small $^1$Center for Combinatorics and LPMC\\
\small Nankai University, Tianjin 300071, China\\
\small Email: lxl@nankai.edu.cn; zhuxy@mail.nankai.edu.cn\\
\small $^2$School of Mathematics and Statistics\\
\small Qinghai Normal University, Xining, Qinghai 810008, China
}

\date{}
\maketitle

\begin{abstract}
A path in an(a) edge(vertex)-colored graph is called a conflict-free path if there exists a color used on only one of its edges(vertices). An(A) edge(vertex)-colored graph is called conflict-free (vertex-)connected if for each pair of distinct vertices, there is a conflict-free path connecting them. For a connected graph $G$, the conflict-free (vertex-)connection number of $G$, denoted by $cfc(G)(\text{or}~vcfc(G))$, is defined as the smallest number of colors that are required to make $G$ conflict-free (vertex-)connected. In this paper, we first give the exact value $cfc(T)$ for any tree $T$ with diameters $2,3$ and $4$. Based on this result, the conflict-free connection number is determined for any graph $G$ with $diam(G)\leq 4$ except for those graphs $G$ with diameter $4$ and $h(G)=2$. In this case, we give some graphs with conflict-free connection number $2$ and $3$, respectively. For the conflict-free vertex-connection number, the exact value $vcfc(G)$ is determined for any graph $G$ with $diam(G)\leq 4$.
\end{abstract}

\section{Introduction}

In this paper, all graphs considered are simple, finite and undirected. We refer to book \cite{BM} for notation and terminology in graph theory not defined here. Among all subjects of graph theory, chromatic theory is no doubt the most arresting. In this paper, we mainly deal with the conflict-free (vertex-) connection coloring of graphs.

In \cite{ELRS}, Even et al. first introduced the hypergraph version of conflict-free (vertex-)coloring. Actually, this coloring emerged as the requirement of the times. It was motivated to solve the problem of assigning frequencies to different base stations in cellular networks. Since then, this coloring has received wide attention due to its practical application value.

Afterwards, Czap et al. introduced the concept of conflict-free connection coloring in \cite{CJV}. In an edge-colored graph, a path is called \emph{conflict-free} if there is at least one color used on exactly one of its edges. This edge-colored graph is said to be \emph{conflict-free connected} if any pair of distinct vertices of the graph are connected by a conflict-free path, and the coloring is called a {\it conflict-free connection coloring}. The \emph{conflict-free connection number} of a connected graph $G$, denoted by $cfc(G)$, is defined as the smallest number of colors required to make $G$ conflict-free connected. There are many results on this topic, for more details, please refer to \cite{CDHJLS,CHL,CJLZ,CJV,DLLMZ}. It is easy to see that $1\leq cfc(G)\leq n-1$ for a connected graph $G$.

Motivated by the above concept, Li et al. \cite{LZZMZJ} introduced the concept of \emph{conflict-free vertex-connection}. A path in a vertex-colored graph is called \emph{conflict-free} if there is a color used on exactly one of its vertices. This vertex-colored graph is said to be \emph{conflict-free vertex-connected} if any two distinct vertices of the graph are connected by a conflict-free path, and the coloring is called a {\it conflict-free vertex-connection coloring}. The \emph{conflict-free vertex-connection number} of a connected graph $G$, denoted by $vcfc(G)$, is defined as the smallest number of colors required to make $G$ conflict-free vertex-connected. In \cite{DS,LZZMZJ,LW}, various results were given in respect of this concept. It has already been obtained that $2\leq vcfc(G)\leq \lceil \log_2(n+1)\rceil$.

We use $S_n$ to denote the \emph{star graph} on $n$ vertices and denote by $T(n_1,n_2)$ the \emph{double star} in which the degrees of its two (adjacent) center vertices are $n_1+1$ and $n_2+1$, respectively. For a connected graph $G$, the \emph{distance} between two vertices $u$ and $v$ is the minimum length of all paths between them, and we write it as $d_G(u,v)$. The \emph{eccentricity} of a vertex $v$ of $G$ is defined by $ecc_G(v)=max_{u\in V(G)}~d_G(u,v)$. The {\it diameter} of $G$ is defined by $diam(G)=max_{v\in V(G)}~ecc_G(v)$ while the {\it radius} of $G$ is defined by $rad(G)=min_{v\in V(G)}~ecc_G(v)$. These parameters have much to do with graph structures and are very significant in the field of graph study. So it stimulates our interest to research on the conflict-free (vertex-)connections of graphs with small diameters.

In this paper, we first give the exact value $cfc(T)$ for any tree $T$ with diameters $2,3$ and $4$. Based of this result, the conflict-free connection number is determined for any graph $G$ with $diam(G)\leq 4$ except for those graphs $G$ with diameter $4$ and $h(G)=2$. In this case, we give some graphs with conflict-free connection numbers $2$ and $3$, respectively. For the conflict-free vertex-connection number, the exact value $vcfc(G)$ is determined for any graph $G$ with $diam(G)\leq 4$.

\section{$cfc$-values for trees with diameters $2,3$ and $4$}\label{tree}

For a connected graph $G$, let $X$ denote the set of cut-edges of $G$, and let $C(G)$ denote the subgraph induced by the set $X$. It is easy to see every component of $C(G)$ is a tree and $C(G)$ is a forest. Let $h(G)=\text{max}~\{cfc(T):T ~\text{is a component of}~ C(G)\}$. In \cite{CJV} and \cite{CHL}, the authors showed the following result.
\begin{lem}\cite{CJV}\label{cfcbound}
For a connected graph $G$, we have $h(G)\leq cfc(G)\leq h(G) + 1.$
Moreover, the bounds are sharp.
\end{lem}
So, $h(G)$ is a crucial parameter to determine the conflict-free connection number of a connected graph $G$. Nevertheless, from the definition of $h(G)$,
determining the value of $h(G)$ depends on determining the conflict-free connection numbers of trees. Therefore, in this section we first give the exact values of the
conflict-free connection numbers of trees with diameters $2,3$ and $4$.
\begin{thm}
For a tree $T$ with diameter $2$ or $3$, we have $cfc(T)=\Delta(T)$.
\end{thm}
\pf It is easy to see that $T$ is a star $S_n$ if and only if it has diameter $2$, and a double star $T(n_1,n_2)(n_1\geq n_2)$ if and only if it has diameter $3$. For the former case, any two edges of $T$ must be colored differently in any conflict-free connection coloring, and thus $cfc(T)=\Delta(T)$. While in the latter case, we can obtain that $cfc(T)=n_1+1=\Delta(T)$ by a similar analysis.\qed

\vspace{0.5cm}
For a tree $T$ of diameter $4$, we denote by $u$ the unique vertex with eccentricity two. The neighbors of $u$ are pendent vertices $w_1,w_2,\cdots,w_\ell$ and $v_1,v_2,\cdots,v_k$ with degrees $p_1\geq p_2\geq\cdots\geq p_k$. Certainly, $k+\ell=d(u)$. In every conflict-free connection coloring $c$ of $T$, the incident edges of every vertex must receive different colors \ding{172}. Without loss of generality, set $c(uv_i)=i(1\leq i\leq k)$ and $c(uw_j)=k+j(1\leq j\leq \ell)$. Observe that if one incident edge of $v_i$ is assigned with color $j$, then color $i$ can not appear on any edge incident with $v_j(1\leq i,j\leq k)$ \ding{173}. Actually, we are seeking for the minimum number of colors satisfying \ding{172} and \ding{173}.

Next we define a vector class $S_r(r\in N^{+})$. We say that an $r$-tuple $(s_1,s_2,s_3,\cdots,s_r)(s_i(1\leq i\leq r)\in N)$ belongs to $S_r$ if we can find a sequence of distinct $2$-tuples $(1,i_{1,1}),(1,i_{1,2}),\cdots,(1,i_{1,s_1}),(2,i_{2,1}),\cdots,(2,i_{2,s_2}), \cdots,$\\ $(r,i_{r,s_r})$ the components of which are all from $[r]$ satisfying that: (1) the two components of every $2$-tuple are different, (2) $(i,j)$ and $(j,i)(1\leq i,j\leq r)$ can not both appear. Note that if $(s_1,s_2,s_3,\cdots,s_r)\in S_r$ then any permutation of its components also belongs to $S_r$. Thus we may suppose $s_1\geq s_2\geq s_3\geq\cdots\geq s_r$.

\begin{lem}\label{tuple}
An $r$-tuple $(s_1,s_2,s_3,\cdots,s_r)(s_i(1\leq i\leq r)\in N)$ belongs to $S_r$ if and only if $\sum\limits_{i=1}^{j}s_i\leq \frac{(2r-1-j)j}{2}(1\leq j\leq r)$.
\end{lem}

\pf First we show the necessity. If $(s_1,s_2,s_3,\cdots,s_r)\in S_r$, then accordingly there is a sequence of $2$-tuples for them according to the definition. Suppose both $(i,j)$ and $(j,i)(1\leq i,j\leq r,i\neq j)$ do not appear. Then, randomly add one of them to the sequence. Repeat this operation until nothing can be added. Finally there are $\frac{(r-1)r}{2}$ $2$-tuples and the corresponding $r$-tuple is $(s'_1,s'_2,s'_3,\cdots,s'_r)$. Assume, to the contrary, there exists a $j$ such that $\sum_{i=1}^{j}s_i> \frac{(2r-1-j)j}{2}(1\leq j\leq r)$. Then $\sum_{i=1}^{j}s'_i\geq \sum_{i=1}^{j}s_i> \frac{(2r-1-j)j}{2}$. Obviously, $j\neq r$. Besides, we have $\sum_{i=j+1}^{r}s'_i\geq \frac{(r-j)(r-j-1)}{2}$ simply by checking the sequence. However, this implies that $\frac{(r-1)r}{2}=\sum_{i=1}^{r}s'_i=\sum_{i=1}^{j}s'_i+\sum_{i=j+1}^{r}s'_i>\frac{(2r-1-j)j}{2}+\frac{(r-j)(r-j-1)}{2}=\frac{(r-1)r}{2}$, a contradiction. Thus the necessity holds.

For the sufficiency, we prove it by applying induction on $r$. When $r=0,1,2$, it is easy to check that if $\sum_{i=1}^{j}s_i\leq \frac{(2r-1-j)j}{2}(1\leq j\leq r)$ for $(s_1,s_2,s_3,$ $\cdots, s_r)$, then this $r$-tuple belongs to $S_r$. Assume that the sufficiency holds for $r=p$. Consider the case $r=p+1$. For $(s_1,s_2,s_3,\cdots,s_{p+1})$, suppose $s_1=p-q$. We distinguish two cases to clarify.

\textbf{Case 1.} $s_{q+1}>s_{q+2}$. In this case, we prove that $(s_2-1,s_3-1,\cdots,s_{q+1}-1, s_{q+2},\cdots,s_{p+1})\in S_p$. When $2\leq j\leq q+1$, we have $\sum_{i=2}^{j}(s_i-1)\leq (j-1)(p-q-1)< \frac{(2p-j)(j-1)}{2}$. When $q+2\leq j\leq p+1$, $\sum_{i=2}^{q+1}(s_i-1)+\sum_{i=q+2}^{j}s_i=\sum_{i=1}^{j}s_i-p\leq \frac{(2p-j+1)j}{2}-p=\frac{(2p-j)(j-1)}{2}$. Therefore, $(s_2-1,s_3-1,\cdots,s_{q+1}-1, s_{q+2},\cdots,s_{p+1})\in S_p$, and so there exists a sequence for it. By adding $(1,p+1),(2,p+1),\cdots,(q,p+1),(p+1,q+1),(p+1,q+2),\cdots,(p+1,p)$ to this sequence, we get a sequence satisfying (1), (2) for $(s_2,\cdots,s_{p+1},s_1)$, implying that $(s_1,s_2,\cdots,s_{p+1})$ belongs to $S_{p+1}$.

\textbf{Case 2.} $s_{q+1}=s_{q+2}$. Suppose $s_r>s_{r+1}=s_{r+2}=\cdots=s_{q+1}=s_{q+2}=\cdots=s_t>s_{t+1}$. Again, we prove that $s'=(s'_1,s'_2,\cdots,s'_p)=(s_2-1,s_3-1,\cdots,s_r-1,s_{r+1},\cdots,s_{t-q+r-1},s_{t-q+r}-1,\cdots,s_t-1,s_{t+1},\cdots,s_{p+1})\in S_p$. Similar to the discussion in \textbf{Case 1}, for $1\leq j\leq r-1$ or $t-1\leq j\leq p$, $\sum_{i=1}^{j}s'_i\leq \frac{(2p-1-j)j}{2}$. Thus if $s'\notin S_p$, the first $j$ such that $\sum_{i=1}^{j-1}s'_i\leq \frac{(2p-j)(j-1)}{2}$ and $\sum_{i=1}^{j}s'_i> \frac{(2p-j-1)j}{2}$ must appear between $r$ and $t-2$. Then we also deduce that $s'_j>p-j$, so $s'_i\geq p-j(j+1\leq i\leq t-1)$. However, this leads to $\frac{(2(p+1)-1-t)t}{2}\geq \sum_{i=1}^{t}s_i=\sum_{i=1}^{t-1}s'_i+p=\sum_{i=1}^{j}s'_i+\sum_{i=j+1}^{t-1}s'_i+p>\frac{(2p-j-1)j}{2}+(t-1-j)(p-j)+p>\frac{(2(p+1)-1-t)t}{2}$,
a contradiction. Thus $s'\in S_p$. By a similar analysis as in \textbf{Case 1}, we can check that $(s_1,s_2,\cdots,s_{p+1})\in S_{p+1}$. The proof is thus complete.\qed

We call the colors from $[k]$ the old colors. In any conflict-free connection coloring of $T$, we denote by $h_i(1\leq i\leq k)$ the number of old colors used on the edges incident with $v_i$ except $uv_i$. Obviously $(h_1,h_2,\cdots,h_k)\in S_k$. In order to add new colors as few as possible, we are actually seeking for the number $a=max(min\{max\{p_i-1-h_i:1\leq i\leq k\}:(h_1,h_2,\cdots,h_k)\in S_k\},0)$.

Let $c_i=p_i-k+i-1(1\leq i\leq k)$, $b=max\{\lceil max\{\sum_{i=1}^{j}\frac{c_i}{j}:1\leq j\leq k\}\rceil,0\}$. Suppose that $max\{\sum_{i=1}^{j}\frac{c_i}{j}:1\leq j\leq k\}$ is obtained by $j=t$. Assume $a<b$. Then $a<\sum_{i=1}^{t}\frac{c_i}{t}$. Thus there exists a $k$-tuple $(h_1,h_2,\cdots,h_k)\in S_k$ such that $h_i\geq p_i-1-a>p_i-1-\sum_{i=1}^{t}\frac{c_i}{t}$. However, this implies that $\sum_{i=1}^{t}h_i>\sum_{i=1}^{t}(p_i-1)-\sum_{i=1}^{t}c_i=\sum_{i=1}^{t}(k-i)=\frac{(2k-1-t)t}{2}$, a contradiction to $(h_1,h_2,\cdots,h_k)\in S_k$ by Lemma \ref{tuple}. Thus $a\geq b$. Next, we only need to construct $(h_1,h_2,\cdots,h_k)\in S_k$ with $b=max\{p_i-1-h_i:1\leq i\leq k\}$. Let $h_i=max\{p_i-1-b,0\}$, it can be easily verified that $(h_1,h_2,\cdots,h_k)$ satisfies our demand. As a result, $a=b$.

Combining Lemma \ref{tuple} with the above analysis, we get the following result.

\begin{thm}
Let $T$ be a tree with diameter $4$, and denote by $u$ its unique vertex with eccentricity two. The neighbors of $u$ are pendent vertices $w_1,w_2,\cdots,w_\ell$ and $v_1,v_2,\cdots,v_k$ with degrees $p_1\geq p_2\geq\cdots\geq p_k$. Then $cfc(T)=max\{k+b,d(u)\}$ where $b=max\{\lceil max\{\sum_{i=1}^{j}\frac{c_i}{j}:1\leq j\leq k\}\rceil,0\}$ and $c_i=p_i-k+i-1(1\leq i\leq k)$.
\end{thm}

\section{Results for graphs with diameters $2,3$ and $4$}

Based on the results in the above section for trees with diameters 2,3, and 4, we are now ready to determine the
$cfc(G)$ and $vcfc(G)$ for graphs with diameters $2,3$ and $4$. At first, we present some auxiliary lemmas that will be used in the sequel.

\begin{lem}\cite{LZZMZJ}\label{vcfc=2}
For a connected graph $G$ of order at least $3$, we have that $vcfc(G)=2$
if and only if $G$ is $2$-connected or $G$ has only one cut-vertex.
\end{lem}

\begin{lem}\cite{LZZMZJ}
For a connected graph $G$, we have $vcfc(G)\leq rad(G)+1$.
\end{lem}

For the conflict-free connection of graphs, the following results have already been obtained.

\begin{lem}\cite{CJV}
For a noncomplete $2$-connected graph $G$, we have $cfc(G)=2$.
\end{lem}

\begin{lem}\cite{CHL}\label{2-edge-connected}
For a noncomplete $2$-edge-connected graph $G$, we have $cfc(G)=2$.
\end{lem}

\begin{lem}\cite{CJV}\label{order2}
If $G$ is a connected graph and $C(G)$ is a linear forest whose each component has an order $2$, then $cfc(G)=2$.
\end{lem}

\begin{lem}\cite{CHL}\label{unique}
Let $G$ be a connected graph with $h(G)\geq 2$. If there exists a unique component $T$ of $C(G)$ such that $cfc(T)=h(G)$, then $cfc(G)=h(G)$.
\end{lem}

\begin{rem}
We have calculated the exact value $cfc(T)$ for any tree $T$ with $diam(T)\leq 4$ in \textbf{Section \ref{tree}}. If $G$ is a connected graph with $diam(G)\leq 4$,
then any component of $C(G)$ must be a tree with diameter no more than four. Thus we can calculate $h(G)$ according to the theorems in \textbf{Section \ref{tree}}.
\end{rem}

For graphs with diameter $2$, we have the following result.

\begin{thm}\label{diameter2}
For a connected graph $G$ with diameter $2$, we have $vcfc(G)=2$ and $cfc(G)=\text{max}~\{2,h(G)\}$.
\end{thm}

\pf Since $G$ has diameter $2$, it is easy to find that $G$ has at most one cut-vertex. According to Lemma \ref{vcfc=2}, $vcfc(G)=2$. If $G$ is $2$-edge-connected, then $cfc(G)=2$ by Lemma \ref{2-edge-connected}. Otherwise, $C(G)$ must be a star, and thus $cfc(G)=\text{max}~\{2,h(G)\}$ by Lemmas \ref{order2} and \ref{unique}.\qed

For graphs with diameter $3$, we have the following result. Recall that a vertex in a block of a graph $G$ is called an \emph{internal vertex} if it is not a cut-vertex of $G$.

\begin{thm}\label{diameter3}
For a connected graph $G$ with diameter $3$, we have that $vcfc(G)\leq3$ and $cfc(G)=max\{2,h(G)\}$ except for the graph depicted in \textbf{Figure 1} which has conflict-free connection number $h(G)+1=3$.
\end{thm}
\begin{figure}[htbp]
\begin{center}
\includegraphics[scale=1.0]{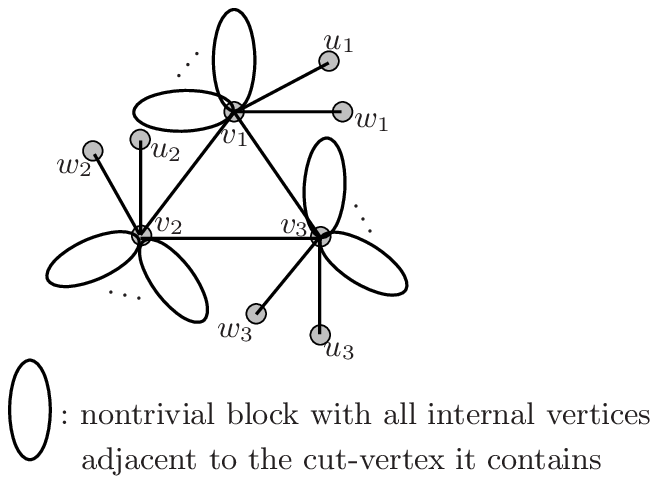}
\caption{}
\end{center}
\end{figure}
\pf Removing all internal vertices of end blocks of $G$, it is easy to check that at most one block is left. Since otherwise if there are two blocks $B_1,B_2$, we can always find two other blocks $C_1,C_2$ such that $V(B_i)\bigcap V(C_i)\neq \phi$ and $V(B_i)\bigcap V(C_{3-i})= \phi(i=1,2)$ and for any two internal vertices $u\in V(C_1),v\in V(C_2)$, every $u$-$v$ path is a $u$-$C_1$-$B_1$-$B_2$-$C_2$-$v$ path. However, this implies that the distance between $u$ and $v$ is at least four, contradicting the fact that $diam(G)=3$. Suppose that $G$ contains no more than one cut-vertex. Then $vcfc(G)=2$ according to Lemma \ref{vcfc=2}. Otherwise for the left block $B_1$, it is bound to contain all cut-vertices of $G$, and we assign the color $3$ to one of them and the color $2$ to all remaining vertices of $V(B_1)$. Other unmentioned vertices of $G$ share the color $1$. It is easy to check that $G$ is conflict-free vertex-connected under this coloring. As a result, $vcfc(G)\leq 3$.

The conflict-free connection number of $G$ has been determined by Lemmas \ref{2-edge-connected}, \ref{order2} and \ref{unique} when $h(G)\leq 1$ or $h(G)\geq 2$ and there exists a unique component $T$ such that $cfc(T)=h(G)$. Thus we only need to consider the remaining cases. This implies that $B_1$ exists and it is nontrivial. Besides, every component of $C(G)$ is a star with its center attached to $B_1$.

Let $h(G)=k$. If $k\geq3$, since $cfc(G)\geq k$, to prove $cfc(G)=k$, we only need to provide a conflict-free connection $k$-coloring of $G$. For each component of $C(G)$, give it a conflict-free connection coloring from $[k]$. As for each nontrivial block, give two of its edges the colors $2$ and $3$ respectively and all others the color $1$. It can be verified that $G$ is conflict-free connected in this way.

When $k=2$, we denote by $n_1$ the number of vertices of $B_1$ and $\ell$ the number of components of $C(G)$. If $\ell<n_1$, then there exists a vertex $v$ of $B_1$ not attached by any component of $C(G)$. Note that since $diam(G)=3$, the subgraph induced by the vertices each of which is attached by some component of $C(G)$ is complete. We only need to give a conflict-free connection $2$-coloring of $G$: The edges of each component of $C(G)$ receive different colors from $[2]$. Randomly choose an edge $e$ of $B_1$ incident with $v$ and each edge for each of other nontrivial blocks, then assign to them the color $2$. The remaining edges are given the color $1$. The checking process is omitted.

For the case $\ell=n_1$, certainly $B_1$ is complete with vertices $v_1,v_2,\cdots,v_{n_1}$. Since $diam(G)=3$, for any end block of $G$, all its internal vertices are adjacent to the cut-vertex it contains. If $n_1\geq4$, we offer a conflict-free connection $2$-coloring of $G$: Assign different colors to the edges of each component of $C(G)$ from $[2]$; give color $2$ to all edges of the path $v_1v_2\cdots v_{n_1}$ and color $1$ to the remaining edges of $B_1$. Observe that each edge of $B_1$ with color $i(i\in[2])$ is contained in a triangle the other two edges of which receive distinct colors. Then pick one edge for each end block and give it color $2$. Other edges are given color $1$. The verification is similar.

Suppose $n_1=3$ with at least one component of $C(G)$ being $P_2$. Choose one of such and give its edge the color $1$. Without loss of generality, assume that this edge is incident with $v_1\in V(B_1)$. Pick one edge of $B_1$ incident with $v_1$ and give it the color $2$, again, other edges of $B_1$ share the same color $1$. We color the edges of other components of $C(G)$ and nontrivial blocks the way as we did in the last paragraph. Obviously, this is a conflict-free connection $2$-coloring for $G$.

If $n_1=3$ and any component of $C(G)$ is a $P_3$, we show that two colors are not enough. Note that there are always two adjacent edges of $B_1$ sharing the same color if only two colors are used. Without loss of generality, suppose that the edges $v_3v_1,v_3v_2$ both have color $1$. Let $v_1u_1,v_2u_2$ have color $1$ and $v_1w_1,v_2w_2$ have color $2$ where these edges are all cut-edges. It is easy to check that there is no conflict-free path between $u_1$ and $u_2$ or $w_1$ and $w_2$ no matter what color the edge $v_1v_2$ is assigned, a contradiction. Thus according to Lemma \ref{cfcbound}, $cfc(G)=h(G)+1=3$.\qed

Finally, we study the conflict-free (vertex-)connection number of graphs with diameter $4$ in the next two results.

\begin{thm}\label{diameter4}
For a connected graph $G$ with diameter $4$, we have that $vcfc(G)\leq 3$, and $cfc(G)=2$ if $h(G)\leq 1$; $cfc(G)=h(G)$ if $h(G)\geq3$.
\end{thm}
\pf Since $G$ has diameter $4$, then after removing all internal vertices of end blocks, the resulting graph has at most one cut-vertex. If there is none, we can give colors as we did in the proof of Theorem \ref{diameter3}. Otherwise, give color $3$ to this cut-vertex $v_1$ and color $2$ to all vertices of blocks incident with $v_1$ except for $v_1$. Finally, assign color $1$ to all remaining vertices. Surely, $G$ is conflict-free vertex-connected under this coloring.

Let $h(G)=k$. If $k\leq 1$, the result follows from Lemmas \ref{2-edge-connected} and \ref{order2}. If $k\geq 3$, we assign to $E(G)$ $k$ colors as we did in the third paragraph of the proof of Theorem \ref{diameter3}. For every pair of distinct vertices $u,v\in V(G)$, any path between them contains the same set of cut-edges. If they belong to the same component of $C(G)$, the conflict-free path is clear. Otherwise, since $diam(G)=4$, there are at most three cut-edges on the path. Thus at least one color of $2$ and $3$ (say $2$) appears at most once. If it does not appear, then we can choose a $u$-$v$ path passing the $2$-colored edge of a nontrivial block and evading all other such edges of the nontrivial blocks it goes through. Else, the desired path is one avoiding all $2$-colored edges of the nontrivial blocks it passes. Thus, $k$ colors are enough in this case.\qed

\begin{cor}
For a connected graph with $diam(G)\leq4$, we have that $vcfc(G)=3$ if and only if $G$ has more than one cut-vertex.
\end{cor}
\pf The result is an immediate corollary of Lemma \ref{vcfc=2}, Theorems \ref{diameter2}, \ref{diameter3} and \ref{diameter4}.\qed

\begin{rem}
If $k=2$, according to Lemma \ref{cfcbound}, we have $2\leq cfc(G)\leq3$. The situation in this case is complicated. Suppose there are exactly $\ell$ components of $C(G)$ with conflict-free connection number 2. Then for each $\ell\geq2$, we give some graphs of diameter $4$ with conflict-free connection numbers 2 and 3, respectively.
\begin{figure}[htbp]
\begin{center}
\includegraphics[scale=1.2]{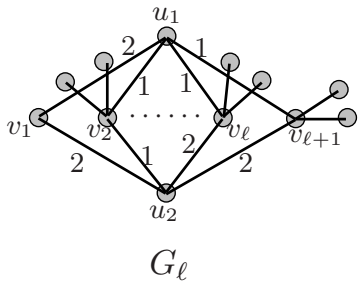}
\caption{: The graph $G_\ell$ with $cfc(G_\ell)=2(\ell\geq2)$.}
\end{center}
\end{figure}

See \textbf{Figure 2} for the graph $G_\ell$ with $cfc(G_\ell)=2(\ell\geq2)$. Each $v_i(2\leq i\leq \ell+1)$ of $G_\ell$ is attached to a $P_3$. We give each such $P_3$ the colors $1$ and $2$ to its two edges, respectively. Besides, give color $1$ to $u_1v_i$ and $2$ to $u_2v_i(3\leq i\leq \ell+1)$. The coloring for other edges are labelled in \textbf{Figure 2}. It is easy to check that this is a conflict-free connection $2$-coloring for $G_\ell$.
\begin{figure}[htbp]
\begin{center}
\includegraphics[scale=1.0]{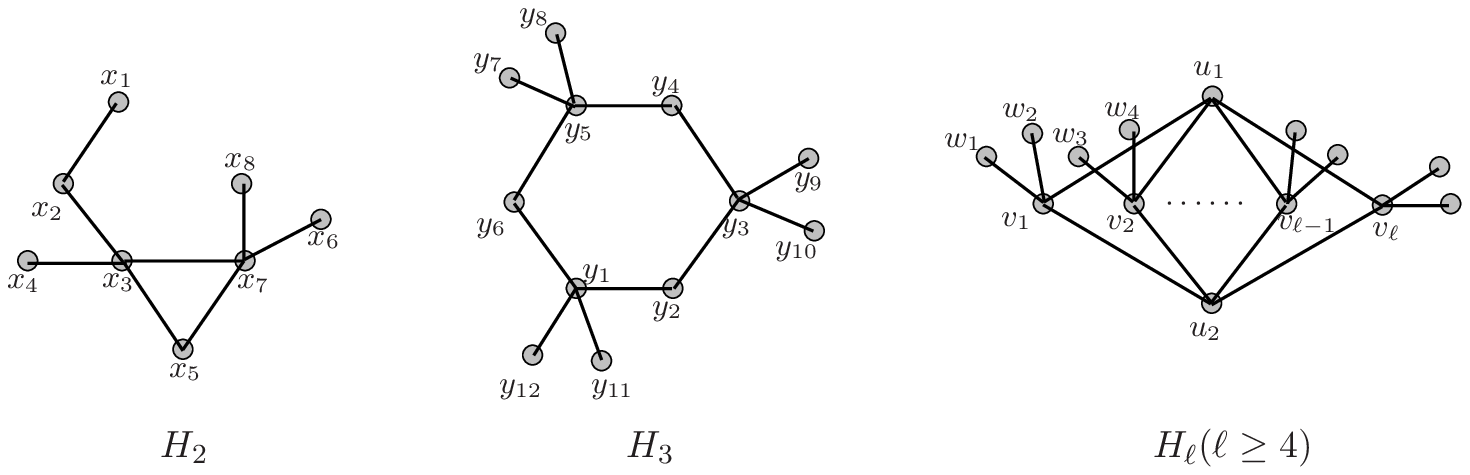}
\caption{: The graph $H_\ell$ with $cfc(H_\ell)=3(\ell\geq2)$.}
\end{center}
\end{figure}

The graph $H_\ell$ with $cfc(H_\ell)=3(\ell\geq2)$ is depicted in \textbf{Figure 3}. Suppose, to the contrary, there exists a conflict-free connection $2$-coloring $c$ for $H_\ell$. When $\ell=2$, without loss of generality, let $c(x_1x_2)=c(x_3x_4)=c(x_6x_7)=1, c(x_2x_3)=c(x_7x_8)=2$. Then if $c(x_3x_7)=1$, to ensure a conflict-free path between $x_4$ and $x_6$, there must be $c(x_3x_5)\neq c(x_5x_7)$. However, there is no conflict-free path between $x_1$ and $x_8$, a contradiction. The case when $c(x_3x_7)=2$ can be dealt with similarly. Thus $cfc(H_2)=3$. With the same method, we can deduce that $cfc(H_3)=3$.

For $H_\ell(\ell\geq4)$, without loss of generality, set $c(v_1w_1)=c(v_2w_3)=1, c(v_1w_2)=c(v_2w_4)=2$. Suppose there exist two monochromatic paths (say $u_1v_1u_2$ and $u_1v_2u_2$) with the same color between $u_1$ and $u_2$. Then there is no conflict-free path between $w_1$ and $w_3$ or $w_2$ and $w_4$, contracting our assumption. If this two monochromatic paths receive different colors, then there is no $w_1$-$w_4$ conflict-free path, a contradiction. Assume that $c(u_1v_1)=c(u_1v_2)\neq c(u_2v_1)=c(u_2v_2)$. For the sake of the existences of conflict-free paths between $w_1$ and $w_3$, $w_2$ and $w_4$, there must be two monochromatic $u_1$-$u_2$ paths with different colors, a contraction to our above analysis. Therefore, $u_1$ and $u_2$ are connected by at most three distinct paths, contradicting with $\ell\geq4$. As a result, $cfc(H_\ell)=3(\ell\geq4)$.
\end{rem}


\begin{thebibliography}{1}

\bibitem{BM}
J. A. Bondy, U. S. R. Murty, Graph Therory, GTM 244, Springer-Verlag, New York, 2008.

\bibitem{CDHJLS}
H. Chang, T.D. Doan, Z. Huang, S. Jendrol¡¯, X. Li, I. Schiermeyer, Graphs with conflict-free connection number two, Graphs Combin., in press.
https://doi.org/10.1007/s00373-018-1954-0.

\bibitem{CHL}
H. Chang, Z. Huang, X. Li, Y. Mao, H. Zhao, On conflict-free connection of graphs, Discrete Appl. Math. 255(2019), 167--182.

\bibitem{CJLZ}
H. Chang, M. Ji, X. Li, J. Zhang, Conflict-free connection of trees, J. Comb. Optim., in press. https://doi.org/10.1007/s10878-018-0363-x.

\bibitem{CJV}
J. Czap, S. Jendrol¡¯, J. Valiska, Conflict-free connection of graphs, Discuss. Math. Graph Theory 38(4), 911--920.

\bibitem{DLLMZ}
B. Deng, W. Li, X. Li, Y. Mao, H. Zhao, Conflict-free connection numbers of line graphs, Lecture Notes in Computer Science No.10627, pp.141--151. Proc. COCOA 2017, Shanghai, China.

\bibitem{DS}
T.D. Doan, I. Schiermeyer, Conflict-free vertex connection number at most $3$ and size of graphs, Manuscript 2018.

\bibitem{ELRS}
G. Even, Z. Lotker, D. Ron, S. Smorodinsky, Conflict-free coloring of simple geometic regions with applications to frequency assignment in cellular networks, SIAM J. Comput. 33(2003), 94--136.


\bibitem{LZZMZJ}
X. Li, Y. Zhang, X. Zhu, Y. Mao, H. Zhao, S. Jendrol', Conflict-free vertex-connections of graphs, Discuss. Math. Graph Theory, in press. DOI: 10.7151/dmgt.2116.

\bibitem{LW}
Z. Li, B. Wu, On the maximum value of conflict-free vertex-connection number of graphs, arXiv:1709.01225 [math.CO].

\end{thebibliography}
\end{document}